\def\PREP{PREP}
\def\ISSAC{ISSAC}
\def\form{PREP}

\ifx\form\PREP
\documentclass[12pt]{amsart}
\fi
\ifx\form\ISSAC
\documentclass{sig-alternate}
\fi

\usepackage{amssymb, colordvi,
multirow, graphicx, algorithmic,setspace}

\newcommand{\bem}[1]{{\Blue{\em #1}}}
\newcommand{\comp}[1]{Y_{#1}}

\def \dId {I^{(d)}}
\def \dXd {X^{(d)}}
\def \dYd {Y^{(d)}}
\def \dZd {Z^{(d)}}
\def \dAd {A^{(d)}}
\def \dLd {L^{(d)}}

\def \dDd {D^{(d)}}

\newcommand{\bC}{{\mathbb C}}
\newcommand{\bQ}{{\mathbb Q}}
\newcommand{\bK}{{\mathbb K}}

\newcommand\Spec{\operatorname{Spec}}
\newcommand\codim{\operatorname{codim}}
\newcommand\NID{\Brown{\operatorname{NID}}}
\newcommand\NPD{\Brown{\operatorname{NPD}}}
\newcommand\IMP{\Brown{\operatorname{IMP}}}
\newcommand\VisibleComponents{\Brown{\operatorname{VisibleComponents}}}
\newcommand\IsInComponent{\Brown{\operatorname{IsInComponent}}}
\newcommand\IsComponent{\Brown{\operatorname{IsComponent}}}
\newcommand\StopCriterium{\Brown{\operatorname{StopCriterium}}}

\newcommand\Ass{\operatorname{Ass}}
\newcommand\VAss{\operatorname{VAss}}

\newcommand{\bff}{{\boldsymbol{f}}}

\newcommand{\bfa}{{\boldsymbol{a}}}

\newcommand{\bfdelta}{{\boldsymbol{\delta}}}
\newcommand{\bfp}{{\boldsymbol{\partial}}}

\newcommand{\xx}{{\boldsymbol{x}}}
\newcommand{\x}{{\mathbf{x}}}
\newcommand{\y}{{\mathbf{y}}}
\newcommand{\z}{{\mathbf{z}}}
\newcommand{\zero}{{\mathbf{0}}}

\newcommand{\p}{{\partial}}

\newcommand{\corank}{{\operatorname{corank}\,}}
 \def\D0{D_\zero}
\def\Span{\operatorname{Span}}
\newcommand{\fm}{{\mathfrak m}}
\newcommand{\maxC}{{\mathcal M}}
\newcommand{\allC}{{\mathcal C}}
\newcommand{\newC}{{\mathcal N}}
\newcommand{\calN}{{\mathcal N}}

\newcommand{\calC}{{\mathcal C}}
\newcommand{\dHd}[2]{{\mathcal H}_{#1,#2}}

\catcode`\@=11\relax \newdimen\p@renwd \font\tenex=cmex10
\setbox0=\hbox{\tenex B} \p@renwd=\wd0
\def\bbordermatrix#1{\begingroup \m@th
\setbox\z@\vbox{\def\\{\crcr\noalign{\kern2\p@\global\let\cr\endline}}%
    \ialign{$##$\hfil\kern2\p@\kern\p@renwd&\thinspace\hfil$##$\hfil
      &&\quad\hfil$##$\hfil\crcr
      \omit\strut\hfil\crcr\noalign{\kern-\baselineskip}%
      #1\crcr\omit\strut\cr}}%
  \setbox\tw@\vbox{\unvcopy\z@\global\setbox\@ne\lastbox}%
  \setbox\tw@\hbox{\unhbox\@ne\unskip\global\setbox\@ne\lastbox}%
  \setbox\tw@\hbox{$\kern\wd\@ne\kern-\p@renwd\left[\kern-\wd\@ne
    \global\setbox\@ne\vbox{\box\@ne\kern2\p@}%
    \vcenter{\kern-\ht\@ne\unvbox\z@\kern-\baselineskip}\,\right]$}%
  \null\;\vbox{\kern\ht\@ne\box\tw@}\endgroup}
\catcode`\@=12\relax

\ifx\form\PREP
\fi
\numberwithin{equation}{section}
\newtheorem{theorem}{Theorem}
\numberwithin{theorem}{section}
\newtheorem{proposition}[theorem]{Proposition}
\newtheorem{lemma}[theorem]{Lemma}

\newtheorem{ex}[theorem]{Example}
\newtheorem{rem}[theorem]{Remark}
\newtheorem{definition}[theorem]{Definition}
\newtheorem{algorithm}[theorem]{Algorithm}

\newenvironment{example}{\begin{ex}\rm}{\end{ex}}
\newenvironment{remark}{\begin{rem}\rm}{\end{rem}}
\newcounter{FNC}[page]
\def\fauxfootnote#1{{\addtocounter{FNC}{2}$^\fnsymbol{FNC}$%
     \let\thefootnote\relax\footnotetext{$^\fnsymbol{FNC}$#1}}}

\begin{document}

\ifx\form\ISSAC
\conferenceinfo{ISSAC'08,} {July 20--23, 2008, Hagenberg, Austria.}
\CopyrightYear{2008}
\crdata{978-1-59593-904-3/08/07}
\renewcommand\Blue[1]{#1}
\renewcommand\Brown[1]{#1}
\fi

\ifx\form\PREP
\def\publname{\scriptsize \Red{Draft of \today} \def\currentvolume{}
\def\currentissue{} \pagespan{1}{60} \PII{}} \copyrightinfo{}{}
\fi

\title{
  Numerical Primary Decomposition
}

\ifx\form\ISSAC
\author{
\alignauthor 
Anton Leykin \\
\affaddr{Institute for Mathematics and its Applications} \\
\affaddr{University of Minnesota} \\
\affaddr{Minneapolis, MN, USA} \\
\email{leykin@ima.umn.edu}
}
\fi
\ifx\form\PREP
\author{Anton Leykin}
\email{leykin@ima.umn.edu}
\urladdr{http://www.ima.umn.edu/\~{}leykin/}
\fi

\maketitle


\begin{abstract}

Consider an ideal $I \subset R = \bC[x_1,\dots,x_n]$ defining a
complex affine variety $X \subset \bC^n$.
We describe the components associated to $I$
by means of {\em numerical primary decomposition}~(NPD).

The method is based on the construction of {\em deflation ideal}
$I^{(d)}$ that defines the {\em deflated variety}  $\dXd$
in a complex space of higher dimension.
For every embedded component there exists $d$ and an isolated component $\dYd$ of $\dId$
projecting onto $Y$. In turn, $\dYd$ can be discovered by
existing methods for prime decomposition, in particular,
the {\em numerical irreducible decomposition}, applied to $\dXd$.

The concept of NPD gives a full description of the scheme $\Spec(R/I)$
by representing each component with a {\em witness set}. We propose an algorithm to produce
a collection of witness sets that contains a NPD and that
can be used to solve the {\em ideal membership problem} for $I$.
\end{abstract}

%
%

\ifx\form\ISSAC
\begin{category}
{Mathematics of Computing}{G.0}{}.
\end{category}

\begin{terms}
Algorithms, Theory.
\end{terms}

\fi



\section{Introduction}

Throughout the paper we use the following notation. Let $I$ be an ideal in the
polynomial ring $R = \bC[\xx] = \bC[x_1,\dots,x_n]$
generated by polynomials $f_1,\dots,f_N \in R$. The ideal $I$ defines the a variety $X=V(I)$,
which, set-theoretically, is the set of points in $\bC^n$ annihilated by all polynomials in ideal $I$.

The basic description of the affine scheme $\Spec(R/I)$
is given by the set of {\em associated prime ideals}
$\Ass(I)$ consisting of all ideals that appear as annihilators of elements of the $R$-module $R/I$.

There exist symbolic algorithms (see, e.g., \cite{
Eisenbud-Huneke-Vasconcelos:primary,
Gianni-Trager-Zacharias,
Caboara-Conti-Traverso,
Decker-Greuel-Pfister,
Steel:primary-decomposition,
Shimoyama-Yokoyama:primary})
to compute an {\em irredundant primary decomposition} of~$I$,
which is by definition a decomposition
\begin{equation} \label{equ:PD}
I = J_1\cap\cdots\cap J_r
\end{equation}
such that all ideals $J_i$ are primary and their
radicals $I_i = \sqrt{J_i} \in \Ass(I)$ are pairwise distinct.

The set of subvarieties defined by the prime ideals
$$\VAss(I) = \{V(P)\ |\ P\in\Ass(I)\}$$
is called the {\em components} associated to~$I$.
A component is said to be {\em isolated} if it is maximal with respect to inclusion
and {\em embedded} if it is not. The existing methods of {\em numerical algebraic geometry}
can ``see'' only the isolated components.

The concept of {\em numerical primary decomposition} (NPD)
that we propose describes {\em all} components in terms of the (generalized) {\em witness sets}
(see Definition \ref{def: witness set}).
Such a witness set for a component $Y\in\VAss(I)$ describes $Y$
completely: one can sample the component, determine its degree, determine whether a given point belongs to $Y$,
etc. Moreover, with a NPD in hand one can solve
the {\em ideal membership problem}, i.e., given a polynomial $f\in R$ decide if
$f\in I$. This is discussed in Subsection~\ref{subsection IMP}.

\medskip

The idea of our method lies in the construction of the {\em deflation ideal} $\dId$
(see Definition \ref{def: deflation ideal}) that defines the {\em deflated variety} $\dXd = V(\dId)$
of order $d$ in a higher-dimensional space. The latter is a stratified vector bundle over $X$
that comes with the natural projection $\pi_d:\dXd\to X$ onto the base. This construction is related
to a deflation technique for 0-dimensional isolated components \cite{LVZ,LVZ-higher,Lec-deflation-02}.

The main theoretical result, Theorem \ref{thm: visible deflation}, shows that
an embedded component $Y \in \VAss(I)$ becomes {\em visible} for some deflation order $d$,
i.e., there is an isolated component $Z \in \VAss(\dId)$ such that $\pi_d(Z) = Y$.

Based on this theorem, 
we outline a straightforward algorithm, Algorithm \ref{alg:VisibleComponents}, to compute all visible components
up to the given order. The advantage of this algorithm over other known techniques is that the problem
of finding all components is reduced to that of finding only isolated components without performing any saturation steps.

Algorithm \ref{alg:NPD} specializes this general algorithm to a numerical method for prime decomposition,
namely, the {\em numerical irreducible decomposition} (NID) \cite{SVW1}.

\medskip

The area of {\em numerical algebraic geometry} comprises novel approaches
to computational algebraic geometry based on the numerical polynomial homotopy continuation methods.
The recent book \cite{Sommese-Wampler-book-05} may serve as a good introduction to the area.

While these methods involve computations that are approximate, a typical output of the algorithms consists
of not only approximations of exact solutions, which are algebraic numbers, but also exact discrete information
about the input. For instance, the main concept introduced in this paper, a {\em numerical} primary decomposition of
an ideal provides {\em exact} data such as the dimensions and the degrees of the associated components.
Another example is a solution to the ideal membership problem given via NPD.
Although we do not provide a certification procedure for this method, such can be developed in theory, provided
that results of our numerical computation can be refined to an arbitrary precision.

\medskip

One advantage of hybrid numerical homotopy continuation techniques over purely symbolic methods,
such as Gr\"obner bases, is that the former is easily parallelizable and the algorithms for computing
the latter are intrinsically serial.

Due to a very small amount
of data being stored on or transferred between the computational nodes
and {\em embarrassing parallelism},
one can easily achieve linear speedups for numerical homotopy
continuation on a parallel computing system with any architecture.
This comes in contrast with the high interdependency of tasks performed in Gr\"obner bases computation
and the phenomena of intermediate expression swell.
(Although there is no doubt that limited speedups are possible in parallel Gr\"obner computation,
the claims of good scalability are substantiated with experiments on preselected classes
of problems and methods based on non-optimal serial algorithms.
See, e.g., \cite{AttardiTraverso} for critique of such claims.)

We live in the age when the clock speed of processors has stopped growing fast
and the computational capacity of computers increases mainly through building
either multicore or distributed systems. This tendency implies both great present
and even better future for the algorithms of numerical algebraic geometry.

\medskip

The paper is structured as follows. After the introduction (Section~1), we make the main definitions
of the paper -- those of deflated ideal and variety -- in Section~2. Next, we overview
the dual space approach to looking at the structure of a polynomial ideal; Section~3 culminates in
the proof of Theorem \ref{thm: visible deflation}, the main theoretical result of the paper.
Section~4 proposes a new numerical representation of an ideal called {\em numerical primary decomposition} (NPD)
via (generalized) {\em witness sets}. In Section~5, we give a skeleton of an NPD algorithm and examples.
We discuss the future of growing meat on the skeleton in conclusion, Section 6.

\medskip

\ifx\form\ISSAC
{\bf Acknowledgements.} I would like to thank IMA (Minneapolis, USA)
for providing an exceptional working environment during 
the thematic year on Applications of Algebraic Geometry,
as well as RISC and RICAM (Linz, Austria) for hosting me for two months 
during the special semester on Gr\"obner bases (2006).
My thanks also go to the referees for their thoughtful comments.
\fi

\section{Deflation}

A variant of the matrix defined below appears in
the {\em deflation} method for regularizing singular isolated solutions of a polynomial
system described in \cite{LVZ-higher}, as well as  the computation of the multiplicity of an isolated point
\cite{DZ-05}.

\begin{definition} The \bem{deflation matrix of order $d$} of
an ideal $I$ generated by polynomials $f_1,\ldots,f_N\in R$  is a matrix $\dAd_I$
with
\begin{itemize}
\item entries in $R$;
\item rows indexed by $\xx^\alpha f_j$, where
$|\alpha| < d$ and $j=1,2,\ldots,N$;
\item columns indexed by
partial differential operators $\bfp^\beta =
\frac{\partial^{|\beta|}}{\partial x_1^{\beta_1}
   \cdots \partial x_n^{\beta_n}}$,
where $|\beta|\leq d$;
\end{itemize}
and with the entry at row
$\xx^\alpha f_j$ and column $\bfp^\beta$ set to be
\begin{equation}
   \partial^\beta\cdot(\xx^\alpha f_j) =
   \frac{\partial^{|\beta|}(\xx^\alpha f_j)}{\partial \xx^\beta}.
\end{equation}
\end{definition}

If a point $\x$ is an isolated solution, in other words, $\{\x\}$ is
0-dimensional irreducible component,
then its multiplicity -- i.e., the vector space dimension of the ring $R/I$ localized at $\x$ --
equals $\corank \dAd_I(\x)$ for a large enough $d$. This
fact is shown, for example, in~\cite{DZ-05} where a method for computing the multiplicity
using $\dAd_I(\x)$ is introduced.

\begin{remark}\label{rem: deflation matrix no 0}
A variant of the deflation matrix
-- the column of $\dAd_I$ labeled with $\bfp^\zero$ is dropped --
is used in the (higher-order) deflation in~\cite{LVZ-higher} to construct
an augmented system of equations for which the multiplicity of a given isolation solution $\x$
of the original system drops. For computational purposes there is no need to keep the aforementioned column,
however, this paper's definition of deflation matrix makes our argument in Section~\ref{sec: dual space} 
more compact.
\end{remark}

\begin{definition}\label{def: deflation ideal}
Let $I=(f_1,\ldots,f_N)\subset R$ and
let $\bfa =\left(a_{\beta}\right)$, $|\beta|\leq
d$, be a vector of indeterminates.

The ideal generated by $f_1,\ldots,f_N$ and the entries of the
vector $\dAd \bfa^T$ in the ring $\bC[\xx,\bfa]$
is called the \bem{deflation ideal} of $I$ of order $d$ and denoted by~\Blue{$\dId$}.

The \bem{deflated variety} of order $d$ is defined
as $$\dXd = V(\dId) \subset \bC^{B(n,d)},$$
where $B(n,d)= \dim \bC[\xx,\bfa] = n+\binom{n+d-1}{d}$.
\end{definition}

Given an ideal $I$, the variety $\dXd$ is well defined in view of the following proposition,
the proof of which is self-contained.

\begin{lemma} \label{lem: poly in ideal}
For every $g \in I$, the deflated variety $\dYd$ of its hypersurface $Y = V(g)$ contains $\dXd = V(\dId)$.
\end{lemma}
\begin{proof}
Let $Q_\bfa = \sum a_\beta \bfp^\beta \in \bC[\bfa, \bfp]$ be the linear differential operator
corresponding to the vector $\bfa = (a_\beta)$ of indeterminates. It is enough to show that
for every point $\x \in X$ and every vector $\bfa$ in the kernel of $\dAd_{\bff}(\x)$
the expression $Q_\bfa g$ vanishes at $\x$. This would imply that the fiber of $\dYd \to Y$ over $\x$
contains the fiber of $\dXd \to X$ over $\x$, thus proving the statement.

Let us write $g = \sum_{\alpha,1\leq i\leq N} c_{i,\alpha}\xx^\alpha f_i$, where $c_{i,\alpha} \in \bC.$
Given a point $\x \in \bC^n$ we can rewrite this as two sums:
$$
g = \sum_{|\alpha|\geq d,1\leq i\leq N} c'_{i,\alpha}(\xx-\x)^\alpha f_i +
    \sum_{|\alpha|<d,1\leq i\leq N} c'_{i,\alpha}\xx^\alpha f_i.
$$
Now, apply the operator $Q_\bfa$ to $g$. The first sum above vanishes at $\x$,
since the order of vanishing of every summand is at least $d+1$.
The second sum vanishes, since $\bfa \in \ker \dAd_{\bff}(\x)$.
\end{proof}

\begin{proposition}\label{prop: variety independence}
If $I_1\subset I_2\subset R$, then $V(\dId_2) \subset V(\dId_1)$ for all $d$.

In particular, the deflation variety $\dXd$
does not depend on the set of generators $f_1,\dots,f_N$ of $I$ chosen to construct $\dAd_I$.
\end{proposition}
\begin{proof}
This follows from Lemma \ref{lem: poly in ideal}, since for a fixed set of generators $f_1,\dots,f_N$ of $I$,
by definition, $$V\left(\dId\right) = V(I)\, \cap \, V\left((f_1)^{(d)}\right) \,\cap\, \cdots \,\cap\, V\left((f_N)^{(d)}\right).$$
\end{proof}

Moreover, a stronger Proposition \ref{prop: deflation independence} holds;
its proof requires a global argument
less transparent than the local argument of Lemma \ref{lem: poly in ideal} and \ref{prop: variety independence}.

\begin{lemma}\label{lem: derivative of Q}
Let $Q_\bfa = \sum_{|\beta|\leq d} a_\beta \bfp^\beta$ and
let $\chi_i Q_\bfa$ is a formal derivative of $Q_\bfa$ with respect to $\p_i$ for $1\leq i\leq n$.

Then $\chi^\alpha Q_\bfa \cdot f_j \in \dId$ for every $j$ and any $\alpha$,  where the deflation ideal $\dId$ is defined using $f_1,\dots,f_N$.
\end{lemma}
\begin{proof}
In case $|\alpha|\geq d$, the expression $\chi^\alpha Q_\bfa$ is a constant in $\bfp$, hence,
$\chi^\alpha Q_\bfa \cdot f_j \in \dId$, since $f_j \in \dId$.

Assume $|\alpha|<d$. By generalized Leibnitz rule
\begin{equation}\label{equ: Q x^alpha f_j}
Q_\bfa\cdot (\xx^\alpha f_j) = \sum_{\gamma} \frac{1}{\gamma!} \frac{\p^{|\gamma|} \xx^\alpha}{\p \xx^\gamma}\chi^\gamma Q_\bfa \cdot f_j.
\end{equation}
The expression $\frac{\p^{|\gamma|} \xx^\alpha}{\p \xx^\gamma}$ is zero if not $\gamma = \alpha$ or $|\gamma|<|\alpha|$.
Assuming that the claim is established for degrees lower than $|\alpha|$,
it follows that $\chi^\gamma Q_\bfa \cdot f_j\in\dId$ for $|\gamma|<|\alpha|$.
Since the left hand side of (\ref{equ: Q x^alpha f_j}) belongs to $\dId$ by definition, it follows that $\chi^\gamma Q_\bfa \cdot f_j$,
the summand corresponding to $\gamma = \alpha$, belongs to $\dId$ as well.
\end{proof}

\begin{proposition}\label{prop: deflation independence}
The deflation ideal $I^{(d)}$ does not depend on the set of generators $I$ chosen to construct it.
\end{proposition}
\begin{proof}
By Lemma \ref{lem: derivative of Q}, we have  $Q_\bfa\cdot (\xx^\alpha f_j)\in \dId$,
therefore, $Q_\bfa\cdot h \in I^{(d)}$ for every $h\in I$, since $h$ is a linear combination of terms $\xx^\alpha f_j$
and the action of $Q_\bfa$ is linear.

The conclusion of the proposition follows immediately.
\end{proof}

\begin{example}\label{exa:(x^2,xyz)}
Let $I = (x_1^2, x_1x_2x_3) \subset \bC[x_1,x_2,x_3]$. The radical $\sqrt{I} = (x_1)$, hence,
the only isolated component is $V(x_1)$.

Compute the first order ($d=1$) deflation: multiply the deflation matrix
$$
A_I^{(1)}(\x) =\ \  \bbordermatrix{
                & \Blue{id} & \Blue{\p_1} & \Blue{\p_2} & \Blue{\p_3} \\
\Blue{x_1^2}    & x_1^2     &  2x_1       & 0           & 0           \\
\Blue{x_1x_2x_3}& x_1x_2x_3 &  x_2x_3     & x_1x_3      & x_1x_2
}
$$
by the column vector $(a_0, a_1, a_2, a_3)^T$. The deflation ideal then is
$$I^{(1)} = (x_1^2,\, x_1x_2x_3,\, a_1x_1,\, a_1x_2x_3+a_2x_1x_3+a_3x_1x_2).$$
Observe that $$\sqrt{I^{(1)}} = (x_1, a_1x_2x_3) = (x_1, a_1) \cap (x_1, x_2) \cap (x_1, x_3),$$
therefore, projecting isolated components of $X^{(1)}$ to $X$ gives two more (embedded) components:
$V(x_1,x_2)$ and $V(x_1,x_3)$.
\end{example}

\begin{remark}\label{rem: stratified bundle}
Denote by $\pi_d:\dXd\to X$ the natural projection induced by
the projection $\pi_d:\bC^{B(n,d)}\to\bC^n$, which maps $(\x,\bfa)\mapsto\x$.

The deflated variety $\dXd$ is a stratified vector bundle over the variety $X$
with the projection $\pi_d$ onto the base. The strata of $X$ over which $\dXd$ is locally trivial
are constructible algebraic subsets of $X$.
\end{remark}

\begin{example}\label{exa: deflation of plane}
Let $L$ be a plane in\, $\bC^n$ of dimension $n-k$ defined by the vanishing
of $I = (x_{j_i})$ for some $1 \leq j_1\leq\cdots\leq j_k\leq n$.

Then it is not hard to establish that
$$\dId = I + ( \{ a_\beta : |\beta|\leq d \mbox{ and } \beta_{j_i} \neq 0 \mbox{ for all } i\} ).$$
Therefore, $\dLd$ is a plane in $\bC^{B(n,k)}$ of dimension
$n-k + \binom{n-k+d-1}{d}$.

Every plane in $\bC^n$ can be brought to the above form with
an affine change of coordinates. This change would result in a linear transformation
on the derivations $\bfp^\beta$; hence, a deflated variety of order $d$ of a plane of codimension $k$
is always a plane of dimension $n-k + \binom{n-k+d-1}{d}$.
\end{example}

\section{Visible components}\label{sec: dual space} 

In this section we describe briefly the dual space (a.k.a. inverse system) approach to the local description
of a point scheme. While from the computational point of view this approach has not been exploited yet beyond
the 0-dimensional case, it comes handy in the proof of the main theorem of the section and the paper.

\subsection{Dual spaces}

\begin{definition}\label{def: dual space}
For $\x\in\bC^n$, let $\Delta_\x^\beta : R\to\bC$ be a linear functional
defined by
$$\Delta_\x^\beta(f) = (\bfp^\beta\cdot f)(\x) = \frac{\p^{|\beta|}f}{\bfp^\beta} (\x),\ \ \ f\in R.$$

The \bem{dual space} $D_\x[I]$ of ideal $I$ at a point $\x$ is the subspace of the
$\bC$-span of functionals $\Delta_\x^\beta$
consisting of the ones that annihilate all polynomials in $I$.
\end{definition}

The dual space defined as above has been introduced under other names in literature;
notably it is the {\em local inverse system}, the term going back to Macaulay \cite{Macaulay:modular-systems}.

The dual space $D_\x[I]$ has a filtration
$$D^{(0)}_\x[I]\subset D^{(1)}_\x[I] \subset D^{(2)}_\x[I] \subset \dots$$
where $\dDd_\x[I]$ is  the set of functionals of order at most $d$.
For computational purposes, it is convenient to think of $\dDd_\x[I]$ as $\ker \dAd_I(\x)$: indeed, given
a vector $\bfa \in \ker \dAd_I(\x)$, one can  view it naturally as the functional
$\sum_{|\beta|\leq d} a_\beta \Delta_\x^\beta$. The natural isomorphism $\ker \dAd_I(\x) \simeq \dDd_\x[I]$
could be shown then by inspection of the definition of the deflation matrix.
In particular, for the case of an isolated point $\x$,
the whole dual space $D_\x[I] \simeq \ker \dAd_I(\x)$ for $d\gg 0$, which is shown, for example,
in \cite[Theorem 1]{DZ-05}.

Let $e_i$ be an $n$-vector with one in the $i$-th position and zeros in the rest.
On the space $D_\x[0]$ of all differential functionals at a point $\x$ we define

\begin{enumerate}

\item operators of \bem{integration} $\delta_i: \dDd_\x[0] \to D^{(d+1)}_\x[0]$, $i=1,\ldots,n$,
induced by the multiplication by $\p_i$ in the Definition \ref{def: dual space}, i.e.:
$$\delta_i(\Delta_\x^\beta)(f) = (\p_i \bfp^\beta \cdot f), \ \ \ f\in R,$$
in other words,
$$\delta_i(\Delta_\x^\beta) = \Delta_\x^{\beta+e_i};$$
\item operators of \bem{differentiation} $\chi_i: D^{(d+1)}_\x[0] \to \dDd_\x[0]$, $i=1,\ldots,n$,
induced by the multiplication by $(x_i-\mathtt{x}_i)$ in the Definition \ref{def: dual space}, i.e.:
\begin{align*}
\chi_i(\Delta_\x^\beta)(f) &= ((x_i-\mathtt{x}_i)\bfp^\beta \cdot f)\\
&= ((\bfp^\beta (x_i-\mathtt{x}_i)+ \beta_i \bfp^{\beta-e_i} )\cdot f),
\end{align*}
therefore, after specializing $\xx=\x$,
$$\chi_i(\Delta_\x^\beta) = \beta_i\Delta_\x^{\beta-e_i}.$$

\end{enumerate}

Since $\Delta_\x^\beta = \bfdelta^\beta \Delta_\x^\zero$, for a fixed $\x$ we may look at the space of linear
functionals $D_\x[0]$ as $\bC[\bfdelta] = \bC[\delta_1,\dots,\delta_n]$
with differentiation operators $\chi_i$ acting as formal derivations $\frac{\p}{\p\delta_i}$.
Note that $D_\x[I] \subset \bC[\delta_1,\dots,\delta_n]$ is stable under this action of $\chi_i$
for every ideal $I\subset R$.

In fact, by (naturally) extending the space of linear functionals to
the formal power series $\bC[[\bfdelta]]$ one may obtain the following fact
(see, e.g., \cite[Proposition 2.6]{Mourrain:inverse-systems}).

\begin{proposition}
The ideals of $R$ are in one-to-one correspondence with the vector subspaces of $\bC[[\bfdelta]]$
stable under differentiation and closed in $(\bfdelta)$-adic topology.
\end{proposition}

For our purposes, given an ideal $I$, we need just a local statement regarding
$(R/I)_\x = (R/I)_{\fm_\x}$, the quotient local ring at a point $\x$, where
$\fm_\x = (x_1-\mathtt{x}_1, \dots, x_n-\mathtt{x}_n) \subset R$.

One may show that the subspaces of $\bC[\bfdelta]$ stable under the differentiation are in one-to-one
correspondence with the ideals of the local ring $R_\x$. The following lemma proves this statement
in one direction.

\begin{lemma}\label{lem:duality}
The image of  $f\in R$ is zero
\begin{enumerate}
  \item in $(R/I+\fm_\x^{d+1})_\x$ \ \ iff \ \  $Q\cdot f = 0$ for all $Q \in \dDd_\x[I]$;
  \item in $(R/I)_\x$ \ \  iff \ \ $Q\cdot f = 0$ for all $Q \in D_\x[I]$.
\end{enumerate}
\end{lemma}
\begin{proof}
Statement (2) follows from (1).

To prove (1), notice that $D_\x[I+\fm_\x^{d+1}]$ is exactly $\dDd_\x[I]$. Indeed, every functional of the latter annihilates
not only $I$, but also all elements of $\fm_\x^{d+1}$; on the other hand, a functional of order larger than $d$ does not
kill the entire $\fm_\x^{d+1}$.
\end{proof}

\begin{remark} In case of an isolated point $\x\in X$, it follows from the lemma
that the dimensions of the $\bC$-spaces  $D_\x[I]$
and $(R/I)_\x$ are equal. Stronger statements are available: for example, see \cite[Theorem 3.1]{LVZ-higher}
and \cite[Theorem 3.2]{Mourrain:inverse-systems}. The latter provides a recipe for constructing the $\fm_\x$-primary component
of the ideal $I$ from a basis of $D_\x[I]$.
\end{remark}

\subsection{A visible deflation of a component}

Recall the deflated variety $\dXd = V(\dId) \subset \bC^{B(n,d)}$ and the projection $\pi_d:\dXd\to X$
in the Definition \ref{def: deflation ideal}.

\begin{definition}\label{def: deflation of component}
We define the \bem{deflation} of order $d$ of a \bem{component} $Y\in \VAss(I)$ as
$\dYd = \overline{\pi_d^{-1}Y^\circ} \subset\dXd$,
where $Y^\circ$ is the subset of \bem{generic points},
that is all smooth points that do not belong to other components
that do not contain~$Y$.
\end{definition}

\begin{proposition}\label{prop: irreducibility of lifting}
A deflation of a component is an irreducible subvariety.
\end{proposition}
\begin{proof}
As was mentioned in Remark \ref{rem: stratified bundle} $\dXd$ is a stratified bundle.
The generic locus $Y^\circ$ is a stratum in the stratification; indeed, the bundle is
locally trivial at every generic point. Since $\pi_d^{-1}Y^\circ$ is an open subset of
$\dYd$, the latter is irreducible.
\end{proof}

\begin{definition}\label{def:visible}
We say that $Y\in\VAss(I)$ is \bem{visible at order $d$},
if $\dYd$ is an isolated component of $\dXd$.
\end{definition}

\begin{theorem}\label{thm: visible deflation} Every component is visible at
some order.
\end{theorem}
\begin{proof}
Note that an isolated component $Y$ is visible for any $d\geq 0$,
since $\dYd = \overline{\pi_d^{-1}Y^\circ}$ can not be contained
in any other component in $\VAss(\dId)$.

Let $I = \bigcap_{Z \in \VAss(I)} J_Z$ be an irredundant primary decomposition
of $I$, where $V(J_Z)=Z$ for all $Z$. Fix an embedded component $Y$ with a generic point
$\y\in Y^\circ$ and consider $I_{>Y} = \bigcap_{Y \subsetneq Z \in \VAss(I)} J_Z$.
Then $I_{>Y} \supsetneq I$, moreover, the local quotient ring $(R/I_{>Y})_\y$ is a proper quotient of
$(R/I)_\y$. Hence, according to Lemma \ref{lem:duality} there is a proper inclusion
$D_\y[I]\supsetneq D_\y[I_{>Y}]$.

For any $Z\in\VAss(I)$ containing $Y$ the fiber of its deflation $\pi_d^{-1}(\y)\cap\dZd$
is contained in $\dDd_\y[I_{>Y}]$, since only the fibers over the generic locus $Z^\circ$ define $\dZd$.
Now, there exists $d$ such that
$$\pi_d^{-1}(\y)\cap\dYd = \dDd_\y[I]\supsetneq \dDd_\y[I_{>Y}] \supset \pi_d^{-1}(\y)\cap\dZd$$
that implies immediately that
$\dZd$ does not contain $\dYd$.
\end{proof}

\begin{remark}
If $Y$ is visible at order $d$ then it is visible at any order $d'\geq d$.
Indeed, the component
$$Y^{(d')} = \overline{\pi_{d'}^{-1}Y^\circ} = \overline{\pi_{d'\to d}^{-1}\left(\dYd\right)^\circ},$$
where $\pi_{d'\to d}:X^{(d')}\to \dXd$ is the natural projection,
has to be isolated if $\dYd$ is isolated.

Since there is a finite number of components there is an order,
at which all components are visible.
\end{remark}

\begin{example}\label{exa:(x^2,xy^2z,xyz^2)}
Consider two ideals of $\bC[x_1,x_2,x_3]$:
\begin{align*}
I &= (x_1^2, x_1x_2x_3)\\
J &= (x_1^2, x_1x_2^2x_3, x_1x_2x_3^2).
\end{align*}

The first order deflation of $I$ has been computed in Example \ref{exa:(x^2,xyz)}. For $J$ we have
$$
A_J^{(1)} =\ \  \bbordermatrix{
                & \Blue{id} & \Blue{\p_1} & \Blue{\p_2} & \Blue{\p_3} \\
\Blue{x_1^2}    & x_1^2     &  2x_1       & 0           & 0           \\
\Blue{x_1x_2^2x_3}& x_1x_2^2x_3 &  x_2^2x_3     & 2x_1x_2x_3      & x_1x_2^2\\
\Blue{x_1x_2x_3^2}& x_1x_2x_3^2 &  x_2x_3^2     & x_1x_3^2      & 2x_1x_2x_3
}
$$
The entries of $A_J^{(1)}(a_0, a_1, a_2, a_3)^T$ together with the original generators of $J$ generate
\begin{eqnarray*}
J^{(1)} = (& x_1^2,\ x_1x_2^2x_3,\ x_1x_2x_3^2, a_1x_1, & \\
               & a_1x_2^2x_3 + 2a_2x_1x_2x_3 + a_3x_1x_2^2, & \\
               & a_1x_2x_3^2 + a_2x_1x_3^2 + 2a_3x_1x_2x_3 & ).
\end{eqnarray*}
It is easy to see that $$\sqrt{J^{(1)}} = \sqrt{I^{(1)}} = (x_1, a_1x_2x_3) = (x_1, a_1) \cap (x_1, x_2) \cap (x_1, x_3),$$
therefore, we need a higher order deflation to distinguish $J$ from $I$.

One may check that the second order deflation uncovers another embedded component in $\VAss(J)$, the origin,
which is not a component in $\VAss(J)$. However, still $\sqrt{J^{(2)}}=\sqrt{I^{(2)}}$, for $I$ the origin is a pseudo-component!

It is not until the third deflation that we see the difference and here is why.
The difference may be seen by comparing the dual spaces at the origin: $D_\zero[I]$ and $D_\zero[J]$. First of all,
$D_\zero[I] \subset D_\zero[J]$, since $J\subset I$.
This inclusion is proper as $\Delta_\zero^{(1,1,1)} \in D_\zero[J] \setminus D_\zero[I]$, a functional of order $3$.
On the other hand, $D^{(2)}_\zero[I] = D^{(2)}_\zero[J]$.
\end{example}

Here we would like to propose the most general algorithm for primary decomposition,
which take the order of deflation $d$ as a parameter and returns all components for $d\gg 0$.

\begin{algorithm} \label{alg:VisibleComponents} $\calN = \VisibleComponents(I,d)$
\begin{algorithmic}
\REQUIRE $I$, ideal of $R={\bK}[\xx]$ where ${\bK}$ is a field of characteristic 0; $d > 0$.
\ENSURE $\calC$, the set components visible at order $d$.
\smallskip \hrule \smallskip
\STATE $\calC = \pi_d\left\{ \mbox{isolated components of }\dId \right\}$.
\smallskip \hrule \smallskip
\end{algorithmic}
\end{algorithm}

The correctness of the algorithm is assured by Theorem \ref{thm: visible deflation};
although we use $\bK = \bC$,
the argument holds for an arbitrary algeraically closed field of characteristic 0. 

In fact, by passing to the algebraic closure, one can make the same conclusion for 
an arbitrary $\bK$ of characteristic 0. The only serious issue in this case is that 
a component that is irreducible over $\bK$ may become reducible over its algebraic closure. 
However, the deflation procedure preserves irreducibility and the property of algebraic closedness 
is exploited only locally in the proof of Theorem \ref{thm: visible deflation}.

\smallskip

While any routine for prime decomposition can be used to compute the isolated components of the deflated variety,
in what follows we concentrate on a numerical approach, which is applicable only for $\bK = \bC$.

\section{Witness sets and numerical primary decomposition}

\subsection{Generalized and classical witness sets}
All numerical algorithms based on homotopy continuation boil down to the computation
of approximations to points of a 0-dimensional variety. That is why for every component $Y$
we need to invent a presentation that would consist of a finite number of points and, perhaps,
some additional (finite) information.

\begin{definition} \label{def: witness set}
A \bem{witness set} $W = W_Y$  of a component $Y\in \VAss(I)$
is a triple $(d, L, w) = (d_Y, L_Y, w_Y)$
consisting of
\begin{enumerate}
\item an order $d$, such that $Y$ is visible at order $d$;

\item a generic  $(\codim \dYd)$-plane $L\subset \bC^{B(n,d)}$;

\item the (finite) set of witness points $w = \dYd \cap L$;
\end{enumerate}
\end{definition}

All items can be presented with finite data: in particular, $L$ can be represented by a linear basis.
We do not include as elements of the witness set generators of $I$, we assume that those are fixed and available.

We also assume that there is a procedure $\dHd{L}{L'}$
that for another generic  $(\codim \dYd)$-plane $L'\subset \bC^{B(n,d)}$ takes the witness points $w_Y$ as input
and produces a new set of witness points $w'_Y$ forming a witness set $(d,L',w'_Y)$. In numerical algebraic geometry
such a procedure is provided by a sufficiently randomized homotopy continuation that deforms $L$ into $L'$
without encountering an intermediate plane that is singular with respect to $\dYd$ and a numerical routine that tracks
the paths starting at the witness points $w_Y$.

We would like to remark that, in principle,
it is enough to store only one
witness point, since the rest can be obtained due to the action of the \bem{monodromy group}, which is
transitive on $w_Y$. In practice, this can be done by following
a random homotopy cycle $\dHd{L}{L}$ a finite number of times.
However, it is, of course, more practical to store the whole set $w_Y$.

\begin{remark}
An isolated component is visible at order $d=0$. In this case, the plane $L$ and the set $w$
give what we would call a {\em classical} witness set, a concept which is used, for example,
throughout \cite{Sommese-Wampler-book-05} where a generalization of it is made (page 237)
in relation to the deflation of an isolated (but multiple) component.
\end{remark}

\begin{example}\label{exa: witness set}
The components of $I = (x_3^2,\,x_3(x_2+x_1^2)) \subset \bC[x_1,x_2,x_3]$
are $Z = V(x_3)$ and $Y = V(x_3, x_2+x_1^2)$.

The isolated component $Z$ can be presented by the witness set using
any line $L$ that is not parallel to $Z$ and such that $Z\cap L \notin Y$.

To represent $Y$, we have to look at the first order deflation:
\begin{align*}
  I^{(1)} &= \big( x_3^2,\,x_3(x_2+x_1^2),\, a_3x_3,\\
          &  \ \ \ \ \ 2a_1x_1x_3+a_2x_3+a_3(x_2+x_1^2) \big)\\
          &  \ \ \ \ \ \ \ \ \ \ \ \ \ \ \subset \bC[x_1,x_2,x_3,a_0,a_1,a_2,a_3]\\
  X^{(1)} &= V(x_3,\,a_3(x_2+x_1^2))\subset \bC^7
\end{align*}
The deflated variety $X^{(1)}$ is 5-dimensional; we take the following 5 equations for $L$:
\begin{align*}
x_2 &= -3x_1+2;\\
a_3 &= x_1-3;\\
a_i &= c_ix_1+d_i,\ \ c_i,d_i\in\bC, \ \ i=0,1,2.
\end{align*}
The first two equations together with the second defining equation of $X^{(1)}$ give
$$
(x_1-3)(x_1-2)(x_1-1) = 0.
$$
The set $\pi_1(X^{(1)}\cap L) = \{(3,-7,0),(2,-4,0),(1,-1,0)\}$
contains projections of two subsets of witness points,
\begin{align*}
\pi_1(w_Z) &= \{(3,-7,0)\},\\
\pi_1(w_Y) &= \{(2,-4,0),(1,-1,0)\},
\end{align*}
of the witness sets of the first order $W_Z=(1,L,w_Z)$ and $W_Y=(1,L,w_Y)$, respectively.
\end{example}

\begin{remark}\label{rem: obtaining dual space}
For a witness point $(\y,\bfa) \in w_Y$, the vector $\bfa$ translates into a functional $Q_\bfa \in \dDd_\y[I]$.
For a fixed generic $\y \in Y$ the set $\{ Q_\bfa \ |\ (\y,\bfa) \in \dYd \}$ equals the dual space $\dDd_\y[I]$.
Therefore, in practice, we can compute $\dDd_\y[I]$ from the witness set $(d,L,w)$ of $Y$ by tracking
a homotopy that creates another witness set $(d,L_i,w_i)$, where the plane $L_i$ is a random plane such that
$\pi_d(L)\cap Y \ni \y$. If the procedure is carried out for sufficiently many $L_i$, then
the functionals $\{ Q_\bfa \ |\ (\y,\bfa) \in w_i \}$ span $\dDd_\y[I]$.
\end{remark}

\subsection{Numerical primary decomposition and the ideal membership problem}
\label{subsection IMP}

\begin{definition}
A collection of witness sets is called a \bem{numerical primary decomposition} (NPD) of $I$
if it contains precisely one witness set for each component in $\VAss(I)$.
\end{definition}

NPD contains exhaustive information about the ideal $I$ and, in particular, the scheme $\Spec(R/I)$ due
to the possibility of solving the \bem{ideal membership problem}. This is so, since for every $Y$
the projection $\pi_{d_Y}(\z)$ of a witness point $\z \in w_Y$ gives a generic point of $Y$ and
in the view of the following:

\begin{theorem}\label{thm:IMP via dual}
A polynomial $g\in R$ is contained in the ideal $I$ iff for all $Y\in \VAss(I)$ and
every (any) generic point $y\in Y$,
all functionals in the dual space $D^{(\deg_\y g)}_\y[I]$ annihilate $g$.
\end{theorem}
\begin{proof}
A polynomial $g\in I$ iff its image in $R/I$ is zero or, equivalently, for all $\x\in X=V(I)$
its image in $(R/I)_\x$ vanishes. It suffices to check the latter statement
for one generic point per component for all components.
\end{proof}

It follows that, in practice, we can solve the ideal membership problem by checking the condition in the theorem
at a finite number of points. Therefore, we have the following algorithm.

\begin{algorithm} \label{alg:IMP}

$b = \IMP(g, I, \calN)$

\begin{algorithmic}

\REQUIRE $I$, ideal of $R$ represented by a finite set of generators;
$\calN$, a NPD of $I$.
\ENSURE $b = \mbox{``}g\in I\mbox{''}$, a boolean value.
\smallskip \hrule \smallskip
\STATE Let $d=\deg g$ and $g=\sum_{|\beta|\leq d} c_\beta\xx^\beta$;
\FORALL{$(d',L,w)\in\calN$}
  \STATE Pick $\x \in \pi_d'(w)$;
  \STATE Compute a linear basis $K$ of $\ker \dAd_I(\x)$;
  \IF{ $\ Q\cdot f \neq 0$ (*) for some $Q\in \dDd_\x[I]$ corresponding to an element of $K$}
  \STATE Return {\bf false};
  \ENDIF
\ENDFOR
\STATE Return {\bf true}.
\smallskip \hrule \smallskip
\end{algorithmic}
\end{algorithm}

According to Theorem \ref{thm:IMP via dual}
$f \in I$ iff there is no $x\in X$, for which the condition (*) holds.
In fact, it suffices to check (*) for a set of generic points (one per component).

\begin{remark}
In case the above algorithm is executed {\em numerically}, i.e., only approximations of
points on the components are generated, we would like to point out several practical issues.

First, having the whole NPD it easy to generate other points on any given component
and recheck the condition (*) at as many points as desired, therefore, lowering the probability
of this algorithm returning an incorrect result due to picking a non-generic point.

Second, since the approximations of generic points can be refined to an arbitrary precision,
the condition (*) can be effectively checked by ``zooming in'' on the exact point for which
the computation is carried out. However, a rigorous certification procedure has to be developed in the future
in the spirit of the {\em alpha-test} for an approximate zero of a univariate polynomial.
\end{remark}

Algorithm \ref{alg:IMP} is practical only for polynomials of low degrees
due to the high complexity of construction of the deflation matrix of order $d$ and computing its kernel.
However, we are confident that improvements can be made as this matrix is highly structured and 
the fact that it is enough to check the condition (*) for only one generic element in the kernel.

\section{Algorithm for numerical\\primary decomposition}

In the description of our NPD algorithm, we assume the following subroutines are at our disposal:

\begin{algorithm} \label{alg:NID}

$\maxC = \NID(I)$

\begin{algorithmic}
\REQUIRE $I$, ideal of $R$.
\ENSURE $\maxC$, the set of classical witness sets $(L_Y,w_Y)$ for all isolated components $Y\in \VAss(I)$.
\end{algorithmic}
\smallskip \hrule \smallskip
There are two approaches to NID in the numerical algebraic geometry: the ``top-down''
method described in detail in \cite[Chapter 15]{Sommese-Wampler-book-05} and ``equation-by-equation''
method the philosophy of which is outlined in \cite[\S 16.2]{Sommese-Wampler-book-05}.
\smallskip \hrule \smallskip
\end{algorithm}

\begin{algorithm} \label{alg: IsInComponent} $b = \IsInComponent(\y,W_Z)$
\begin{algorithmic}
\REQUIRE $\y \in \bC^n$;
         $W_Z=(d,L,w)$ a witness set for $Z\in\VAss(I)$.
\ENSURE $b = \mbox{``\,$y \in Z$''}$, a boolean value.
\end{algorithmic}
\smallskip \hrule \smallskip
Pick a generic ($\codim Z$)-plane $M\subset\bC^n$ and a ($\dim L$)-plane $L'\subset \bC^{B(n,d)}$
such that $\pi_d(L') = M \ni \y$.

Use a procedure similar to the usual containment test routine \cite[\S 15.1]{Sommese-Wampler-book-05},
i.e.,  track the points $w$ along a generic homotopy $\dHd{L}{L'}$.

Return ``\,$\y \in \pi_d(\dHd{L}{L'}(w))$''.
\smallskip \hrule \smallskip
\end{algorithm}

There are also two subroutines, for which finding efficient algorithms is an open problem:
\begin{itemize}
  \item $\StopCriterium(d, I, \calN)$ implements a termination criterion that guarantees that all components
  in $\VAss(I)$ are visible at order less than $d$.
  \item $\IsComponent((d,L,w), I, \calN)$ is used to filter out witness sets $(d,L,w)$ that represent \bem{false components}
  that appear due to singularities. As a parameter it takes partial NPD $\calN$ that includes the witness sets for all
  components visible at order $d$ of dimensions higher than the dimension of the alleged component.
\end{itemize}
Both routines will be discussed later in this section.

Now we are ready to outline the main algorithm of this paper.

\begin{algorithm} \label{alg:NPD} $\calN = \NPD(I)$
\begin{algorithmic}
\REQUIRE $I$, ideal of $R$.
\ENSURE $\calN$, the set of witness sets for all $Y\in \VAss(I)$ visible at order $d$.
\smallskip \hrule \smallskip

\STATE $\calN = \emptyset$;
\REPEAT
 \STATE $\begin{array}{l}
    \newC = \big\{ (L,w) \in \NID(\dId) :\\
    \ \ \ \ \ \ \mbox{not }\IsInComponent(\pi_d(w_1),W) \mbox{ for all } W \in \allC \big\};
    \end{array}$
 \FORALL{$(L,w) \in \newC$ in an order

of decreasing $\dim \comp{(L,w)}$}
   \IF{ $\IsComponent((d,L,w), I, \calN)$ }
      \STATE $\calN = \calN \cup \{ (d,L,w) \}$;
   \ENDIF
 \ENDFOR
 \STATE $d = d+1$;
\UNTIL{ $\StopCriterium(d, I, \calN)$ }

\smallskip \hrule \smallskip
\end{algorithmic}
\end{algorithm}

The order of deflation sufficient to discover all components can be bound by the maximum of
the regularities of the (local) Hilbert functions at the points of $X$.
For the latter a crude bound exists, which is doubly exponential in the number of variables.
Obviously, one can not use $\StopCriterium$ based on this bound for practical purposes.

In reality, for many nontrivial examples {\em all}\,
embedded components are discovered by the deflation of order as low as 1. While finding a reasonable termination
criterion is the matter of the future, currently it makes sense to run a {\em truncated} computation with
$$\StopCriterium(d, I, \calN)\ =\  \mbox{``}d>d_{max}\mbox{''},$$
where $d_{max}$ is the maximal deflation order considered.

As to $\IsComponent$, first of all, we remark that elimination of fake components from $\calN$
may also be done at the end of the algorithm.
There is a way to do this by checking whether for large enough $d$ all functionals in the $\dDd_\y[I]$ for a generic $\y$
in an alleged component $Y$ ``come'' from the components that contain $Y$.
We do not describe the procedure here as it is quite technical and not practical at the moment,
since the question how large $d$ should be relates to the question of finding a good $\StopCriterium$.

\begin{remark}\label{rem: unfiltered NPD}
Establishing an efficient $\StopCriterium$  has a higher priority
(over $\IsComponent$), since $\calN$ containing additional witness sets
of fake components can be used instead of a true NPD for many tasks.
In particular, our $\IMP$ routine (Algorithm~\ref{alg:IMP}) would still work.
\end{remark}

Below is an example of a ``truncated computation'',  where only the first order deflation is computed.
\begin{example}
Consider the cyclic 4-roots problem:
\begin{align*}
I &=  \big( x_1+x_2+x_3+x_4,\ x_1x_2+x_2x_3+x_3x_4+x_4x_1,\\
  & x_1x_2x_3+x_2x_3x_4+x_3x_4x_1+x_4x_1x_2,\ x_1x_2x_3x_4-1 \big )\,.
\end{align*}

The calculation of associated primes via symbolic software (we used {\em Macaulay~2}~\cite{M2www}) gives:
$$\begin{array}{rll}
  \Ass(I) =  \big\{ & (x_2+x_4,x_1+x_3,x_3x_4+1),\\
                    & (x_2+x_4,x_1+x_3,x_3x_4-1),\\
                    & (x_4-1,x_3+1,x_2+1,x_1-1),\\
                    & (x_4-1,x_3-1,x_2+1,x_1+1),\\
                    & (x_4+1,x_3+1,x_2-1,x_1-1),\\
                    & (x_4+1,x_3-1,x_2-1,x_1+1),\\
                    & (x_3+x_4,x_2+x_4,x_1-x_4,x_4^2+1),\\
                    & (x_3-x_4,x_2+x_4,x_1+x_4,x_4^2+1)& \big\}
\end{array}$$
The first two ideals correspond to the irreducible curves that are the two isolated components.
The rest are embedded 0-dimensional components; note that the last two ideals are irreducible
over the ground field $\bQ$, but not $\bC$.

Over complex numbers, there are 8 embedded components
that are all visible at order~1;
The numerical computation of the irreducible components of the first deflation
finds all components.
An excerpt from the numerical output is given below: we list the projections of 
the witness points for the components $\VAss(I^{(1)})$. 
\end{example}

{
\begin{verbatim}
>>> projections of witness points for
>>> component #1:
[x1 = -4.4882+2.0260*I, x2 = .18509+.83550e-1*I, 
x3 = 4.4882-2.0260*I, x4 = -.18509-.83550e-1*I]
[x1 = .52885e-1-.87608*I, x2 = -.68654e-1-1.1373*I, 
x3 = -.52885e-1+.87608*I,x4 = .68654e-1+1.1373*I]
[x1 = -.15083+.49191*I, x2 = .56975+1.8582*I, 
x3 = .15083-.49191*I, x4 = -.56975-1.8582*I]
[x1 = .41488+.24720*I, x2 = -1.7788+1.0599*I, 
x3 = -.41488-.24720*I, x4 = 1.7788-1.0599*I]
>>> component #2:
[x1 = -.95775+.36799*I, x2 = -.90980-.34957*I, 
x3 = .95775-.36799*I, x4 = .90980+.34957*I]
[x1 = .71538+.12328*I, x2 = 1.3576-.23395*I, 
x3 = -.71538-.12328*I, x4 = -1.3576+.23395*I]
[x1 = -3.7686+1.7072*I, x2 = -.22017-.99738e-1*I, 
x3 = 3.7686-1.7072*I, x4 = .22017+.99738e-1*I]
[x1 = -.16036-.30943*I, x2 = -1.3202+2.5476*I, 
x3 = .16036+.30943*I, x4 = 1.3202-2.5476*I]
>>> component #3:
[x1 = -1.0-.53734e-17*I, x2 = 1.0-.20045e-16*I, 
x3 = 1.0+.89149e-17*I, x4 = -1.0+.18026e-17*I]
...
>>> component #10:
[x1 = -.59351e-17+1.0*I, x2 = -.46995e-16+1.0*I, 
x3 = .16158e-16-1.0*I, x4 = .22439e-16-1.0*I]
\end{verbatim}
}
Note that there are 2 witness sets of 4 points 
corresponding to the 2 isolated curves and 8 singletons for the embedded points.

See the webpage \cite{NPDwww} for the scripts in {\em Macaulay~2} 
and {\em Maple} (using {\em PHCmaple} package~\cite{PHCMAPLEwww})
that perform prime decomposition and numerical irreducible decomposition, 
respectively, for the first deflation ideal in this example.


\section{Discussion and conclusion}

We consider this paper as one laying a theoretical foundation to the method that at this point works only on small examples.
However, it is our believe that this technique would be able to solve problems unsolvable
by purely symbolic methods in the future. The improvements are expected to be made both in the software and in the theory.

We remark that the software in the area of numerical algebraic geometry is as young as the area itself.
For the purposes of numerical irreducible decomposition there exist only two software options: {\em PHCpack}~\cite{V99}
that we use via {\em PCHmaple}~\cite{PHCMAPLEwww} and {\em Bertini}~\cite{Bertini}.
The practical computation using the ideas in this paper is limited by the capabilities of these software systems;
we expect the implementations of numerical irreducible decomposition algorithms to improve.
Both {\em PHCpack} and {\em Bertini} move towards the throughout parallelization; as we argued in the introduction,
easy parallelization is a crucial feature of numerical methods that distinguishes them from the symbolic ones.

The future theoretical work should, in particular, concentrate on the construction of special homotopy methods
to tackle deflated ideals (systems),
which possess an obvious multihomogeneous structure: they are linear in the additional variables.
Also, while the global algorithms such as NID have been well established, there are still no efficient
local procedures, e.g., for determining the local dimension at a given point on a variety.
The same can be said about the local dual space computation: while it is possible to compute the truncation
at some degree using the deflation matrix,
an efficient description of the whole (possibly infinite-dimensional)
dual space and ways to create such description are yet to be found.



\begin{thebibliography}{10}

\bibitem{AttardiTraverso}
G.~Attardi and C.~Traverso.
\newblock {Strategy-accurate parallel Buchberger algorithms}.
\newblock {\em J. Symbolic Comput.}, 21(4-6):411--425, 1996.

\bibitem{Bertini}
D.~J. Bates, J.~D. Hauenstein, A.~J. Sommese, and C.~W. Wampler.
\newblock Bertini: software for numerical algebraic geometry.
\newblock Available at http://www.nd.edu/$\sim$sommese/bertini.

\bibitem{Caboara-Conti-Traverso}
M.~Caboara, P.~Conti, and C.~Traverso.
\newblock Yet another ideal decomposition algorithm.
\newblock In {\em Applied algebra, algebraic algorithms and error-correcting
  codes (Toulouse, 1997)}, volume 1255 of {\em Lecture Notes in Comput. Sci.},
  pages 39--54. Springer, Berlin, 1997.

\bibitem{DZ-05}
B.~Dayton and Z.~Zeng.
\newblock Computing the multiplicity structure in solving polynomial systems.
\newblock In M.~Kauers, editor, {\em Proceedings of the 2005 International
  Symposium on Symbolic and Algebraic Computation}, pages 116--123. ACM, 2005.

\bibitem{Decker-Greuel-Pfister}
W.~Decker, G.-M. Greuel, and G.~Pfister.
\newblock Primary decomposition: algorithms and comparisons.
\newblock In {\em Algorithmic algebra and number theory (Heidelberg, 1997)},
  pages 187--220. Springer, Berlin, 1999.

\bibitem{Eisenbud-Huneke-Vasconcelos:primary}
D.~Eisenbud, C.~Huneke, and W.~Vasconcelos.
\newblock Direct methods for primary decomposition.
\newblock {\em Invent. Math.}, 110(2):207--235, 1992.

\bibitem{Gianni-Trager-Zacharias}
P.~Gianni, B.~Trager, and G.~Zacharias.
\newblock Gr\"obner bases and primary decomposition of polynomial ideals.
\newblock {\em J. Symbolic Comput.}, 6(2-3):149--167, 1988.
\newblock Computational aspects of commutative algebra.

\bibitem{M2www}
D.~R. Grayson and M.~E. Stillman.
\newblock Macaulay 2, a software system for research in algebraic geometry.
\newblock Available at http://www.math.uiuc.edu/Macaulay2/.

\bibitem{Lec-deflation-02}
G.~Lecerf.
\newblock Quadratic {N}ewton iteration for systems with multiplicity.
\newblock {\em Found. Comput. Math.}, 2:247--293, 2002.

\bibitem{NPDwww}
A.~Leykin.
\newblock {Numerical primary decomposition (webpage)}.
\newblock www.math.umn.edu/$\sim$leykin/NPD.

\bibitem{PHCMAPLEwww}
A.~Leykin.
\newblock {PHCmaple: A Maple interface to the numerical homotopy algorithms in
  PHCpack}.
\newblock \newline www.math.umn.edu/$\sim$leykin/PHCmaple.

\bibitem{LVZ}
A.~Leykin, J.~Verschelde, and A.~Zhao.
\newblock Newton's method with deflation for isolated singularities of
  polynomial systems.
\newblock {\em Theoretical Computer Science}, 359(1-3):111--122, 2006.

\bibitem{LVZ-higher}
A.~Leykin, J.~Verschelde, and A.~Zhao.
\newblock Higher-order deflation for polynomial systems with isolated singular
  solutions.
\newblock In A.~Dickenstein, F.-O. Schreyer, and A.~J. Sommese, editors, {\em
  Algorithms in Algebraic Geometry}, volume 146 of {\em The IMA Volumes in
  Mathematics and its Applications}. Springer, 2008.

\bibitem{Macaulay:modular-systems}
F.~S. Macaulay.
\newblock {\em The algebraic theory of modular systems}.
\newblock Cambridge Mathematical Library. Cambridge University Press,
  Cambridge, 1994.
\newblock Revised reprint of the 1916 original, With an introduction by Paul
  Roberts.

\bibitem{Mourrain:inverse-systems}
B.~Mourrain.
\newblock Isolated points, duality and residues.
\newblock {\em J. Pure Appl. Algebra}, 117/118:469--493, 1997.
\newblock Algorithms for algebra (Eindhoven, 1996).

\bibitem{Shimoyama-Yokoyama:primary}
T.~Shimoyama and K.~Yokoyama.
\newblock Localization and primary decomposition of polynomial ideals.
\newblock {\em J. Symbolic Comput.}, 22(3):247--277, 1996.

\bibitem{SVW1}
A.~Sommese, J.~Verschelde, and C.~Wampler.
\newblock Numerical decomposition of the solution sets of polynomial systems
  into irreducible components.
\newblock {\em SIAM J.\ Numer.\ Anal.}, 38(6):2022--2046, 2001.

\bibitem{Sommese-Wampler-book-05}
A.~J. Sommese and C.~W. Wampler, II.
\newblock {\em The numerical solution of systems of polynomials}.
\newblock World Scientific Publishing Co. Pte. Ltd., Hackensack, NJ, 2005.

\bibitem{Steel:primary-decomposition}
A.~Steel.
\newblock Conquering inseparability: primary decomposition and multivariate
  factorization over algebraic function fields of positive characteristic.
\newblock {\em J. Symbolic Comput.}, 40(3):1053--1075, 2005.

\bibitem{V99}
J.~Verschelde.
\newblock Algorithm 795: {PHC}pack: A general-purpose solver for polynomial
  systems by homotopy continuation.
\newblock {\em ACM Trans. Math. Softw.}, 25(2):251--276, 1999.
\newblock Software available at {\tt http://www.math.uic.edu/{\~{}}jan}.

\end{thebibliography}

\bibliographystyle{abbrv}

\def\cprime{$'$}

\end{document}